\documentclass[a4paper,11pt]{amsart}
\usepackage{graphicx}
\usepackage{mathptmx}
\usepackage{amsmath}
\usepackage{amssymb}
\usepackage{enumitem}
\usepackage{xcolor}
\usepackage{xparse}
\NewDocumentCommand{\eulerian}{omm}
 {%
  \genfrac<>{0pt}{}{#2}{#3}%
  \IfValueT{#1}{_{\!#1}}%
 }

 \newmuskip\pFqmuskip

\newcommand*\pFq[6][8]{%
  \begingroup % only local assignments
  \pFqmuskip=#1mu\relax
  \mathchardef\normalcomma=\mathcode`,
  \mathcode`\,=\string"8000
  \begingroup\lccode`\~=`\,
  \lowercase{\endgroup\let~}\pFqcomma
  {}_{#2}F_{#3}{\left(\genfrac..{0pt}{}{#4}{#5}\bigg|#6\right)}%
  \endgroup
}
\newcommand{\pFqcomma}{{\normalcomma}\mskip\pFqmuskip}

\newtheorem{theorem}{Theorem}

\newtheorem{corollary}[theorem]{Corollary}

\begin{document}

\title[A note on degenerate generalized Laguerre polynomials and Lah numbers]{A note on degenerate generalized Laguerre polynomials and Lah numbers}

\author{Taekyun  Kim}
\address{Department of Mathematics, Kwangwoon University, Seoul 139-701, Republic of Korea}
\email{kwangwoonmath@hanmail.net}

\author{Dmitry V. Dolgy}
\address{Department of Mathematical Methods in Economy, Far Eastern Federal University, 690950 Vladivostok,	Russia}
\email{dvdolgy@gmail.com }

\author{Dae San Kim}
\address{Department of Mathematics, Sogang University, Seoul 121-742, Republic of Korea}
\email{dskim@sogang.ac.kr}

\author{Hye Kyung Kim}
\address{Department Of Mathematics Education, Daegu Catholic University, Gyeongsan 38430, Republic of Korea}
\email{hkkim@cu.ac.kr}

\author{Seong Ho Park}
\address{Department of Mathematics, Kwangwoon University, Seoul 139-701, Republic of Korea}
\email{abcd2938471@kw.ac.kr}


\begin{abstract}
The aim of this paper is to introduce the degenerate generalized Laguerre polynomials as the degenerate version of the generalized Laguerre polynomials and to derive some properties related to those polynomials and Lah numbers, including an explicit expression, a Rodrigues' type formula and expressions for the derivatives. The novelty of the present paper is that it is the first paper on degenerate versions of orthogonal polynomials.
\end{abstract}

\maketitle

\section{Introduction}
The generalized Laguerre polynomials are classical orthogonal polynomials which are orthogonal with respect to the gamma distribution $e^{-x} x^{\alpha}dx$ on the interval $ \left(0,\infty\right)$.
The generalized Laguerre polynomials are widely used in many problems of quantum mechanics, mathematical physics and engineering.
In quantum mechanics, the Schr\"{o}dinger equation for the hydrogen-like atom is exactly solvable by separation of variables in spherical coordinates. The radial part of the wave function is a generalized Laguerre polynomial [14].
In mathematical physics, vibronic transitions in the Franck-Condon approximation can also be described by using Laguerre polynomials [6].
In engineering, the wave equation is solved for the time domain electric field integral equation for arbitrary shaped conducting structures by expressing the transient behaviors in terms of Laguerre polynomials [4].\par
The aim of this paper is to introduce the degenerate generalized Laguerre polynomials as the degenerate version of the generalized Laguerre polynomials and to derive some properties related to those polynomials and Lah numbers. In more detail, we obtain an explicit formula and a Rodrigues' type formula for the degenerate Laguerre polynomials. We also get explicit expressions for the degenerate generalized Laguerre polynomial for $\alpha=-1$, an identity involving Lah numbers, the falling factorial moment of the degenerate Poisson random variable with parameter $\alpha$, and expressions for the derivatives of the degenerate generalized Laguerre polynomials. \par
We should mention here that degenerate versions of many special numbers and polynomials have been explored and many interesting results have obtained in recent years [8,10,12]. Furthermore, these have been done not only for special numbers and polynomials but also for transcendental functions like gamma functions [11]. The novelty of the present paper is that this is the first paper which treats degenerate versions of orthogonal polynomials. For the rest of this section, we will recall some necessary facts that will be used throughout this paper.

The Lagueree polynomial $L_{n}(x)$ satisfies the second-order linear differential equation
\begin{displaymath}
	xy''+(1-x)y'+ny=0,\quad(\mathrm{see}\ [16]),
\end{displaymath}
while the generalized Laguerre polynomial (or the associated Laguerre polynomial) $L_{n}^{(\alpha)}(x)$ satisfies the second-order linear differential equation
\begin{displaymath}
	xy''+(\alpha+1-x)y'+ny=0,\quad(\alpha\in\mathbb{R}).
\end{displaymath}
The Rodrigues'  formula of the Laguerre polynomial $L_{n}(x)$ is given by
\begin{equation}
	L_{n}(x)=\frac{e^{x}}{n!}\frac{d^{n}}{dx^{n}}\big(e^{-x}x^{n}\big)=\frac{1}{n!}\bigg(\frac{d}{dx}-1\bigg)^{n}x^{n},\label{1}
\end{equation}
while that of the generalized Lagurerre polynomial $L_{n}^{(\alpha)}(x)$ is given by
\begin{align}
L_{n}^{(\alpha)}(x)\ &=\ \frac{1}{n!}x^{-\alpha}e^{x}\frac{d^n}{dx^n}\big(e^{-x}x^{n+\alpha}\big) \label{2}\\
&=\ x^{-\alpha}\frac{1}{n!}\bigg(\frac{d}{dx}-1\bigg)^{n}x^{n+\alpha},\quad (\mathrm{see}\ [2,9,16,17]).\nonumber
\end{align}
The generating function of generalized Laguerre polynomials is given by
\begin{equation}
\sum_{n=0}^{\infty}L_{n}^{(\alpha)}(x)t^{n}=\frac{1}{(1-t)^{\alpha+1}}e^{-x\frac{t}{1-t}},\quad(\mathrm{see}\ [9,16,17]).\label{3}	
\end{equation}
From \eqref{3}, we get
\begin{equation}
L_{n}^{(\alpha)}(x)=\sum_{i=0}^{n}(-1)^{i}\binom{n+\alpha}{n-i}\frac{x^{i}}{i!},\quad(\mathrm{see}\ [9,17]).\label{4}
\end{equation}
Note that
\begin{align*}
L_{0}^{(\alpha)}(x)&=1 	\\
L_{1}^{(\alpha)}(x)&=-x+(\alpha+1) \\
L_{2}^{(\alpha)}(x)&=\frac{x^{2}}{2}-(\alpha+2)x+\frac{(\alpha+1)(\alpha+2)}{2}, \\
L_{3}^{(\alpha)}(x)&=-\frac{x^{3}}{6}+\frac{\alpha+3}{2}x^{2}-\frac{(\alpha+2)(\alpha+3)}{2}x+\frac{(\alpha+1)(\alpha+2)(\alpha+3)}{6},\ \dots.
\end{align*}
The rising factorial sequence is defined as
\begin{displaymath}
	\langle x\rangle_{0}=1,\quad \langle x\rangle_{n}=x(x+1)\cdots(x+n-1),\quad (n\ge 1),
\end{displaymath}
while the falling factorial sequence is defined as
\begin{displaymath}
	(x)_{0}=1,\quad (x)_{n}=x(x-1)\cdots(x-n+1),\quad (n\ge 1),\quad (\mathrm{see}\ [1-3,5,7-13,15-18]).
\end{displaymath}
We note that the Lah numbers are defined by
\begin{equation}
\langle x\rangle_{n}=\sum_{k=0}^{n}L(n,k)(x)_{k},\quad (n\ge 0),\quad(\mathrm{see}\ [3,5,8,13,15]).\label{5}
\end{equation}
From \eqref{5}, we can easily derive the following equation.
\begin{equation}
\frac{1}{k!}\bigg(\frac{t}{1-t}\bigg)^{k}=\sum_{n=k}^{\infty}L(n,k)\frac{t^{n}}{n!},\quad(\mathrm{see}\ [3,5,13,17]).\label{6}	
\end{equation}
For any $\lambda\in\mathbb{R}$, the degenerate exponential function is defined by
\begin{equation}
e_{\lambda}^{x}(t)=\sum_{n=0}^{\infty}(x)_{n,\lambda}\frac{t^{n}}{n!},\quad(\mathrm{see}\ [10,11])\label{7}
\end{equation}
where $(x)_{0,\lambda}=1,\ (x)_{n,\lambda}=x(x-\lambda)\cdots(x-(n-1)\lambda)$, $(n\ge 1)$.
For $x=1$, we use the notation by $e_{\lambda}(t)=e_{\lambda}^{1}(t)$. \par

\section{Degenerate generalized Laguerre polynomials}
For any $\alpha\in\mathbb{R}$, we consider the degenerate generalized Laguerre polynomials given by
\begin{equation}
\frac{1}{(1-t)^{\alpha+1}}e_{\lambda}\bigg(-x\frac{t}{1-t}\bigg)=\sum_{n=0}^{\infty}L_{n,\lambda}^{(\alpha)}(x)t^{n},\quad |t|<1.\label{8}	
\end{equation}
From \eqref{7}, we note that
\begin{align}
	&\frac{1}{(1-t)^{\alpha+1}}e_{\lambda}\bigg(-x\frac{t}{1-t}\bigg)=\frac{1}{(1-t)^{\alpha+1}}\sum_{m=0}^{\infty}(1)_{m,\lambda}(-1)^{m}x^{m}\frac{1}{m!}\bigg(\frac{t}{1-t}\bigg)^{m}\label{9} \\
	&\quad =\sum_{m=0}^{\infty}(1)_{m,\lambda}(-1)^{m}x^{m}\frac{1}{m!}t^{m}\bigg(\frac{1}{1-t}\bigg)^{m+\alpha+1}\nonumber \\
	&\quad = \sum_{m=0}^{\infty}(1)_{m,\lambda}(-1)^{m}x^{m}\frac{t^{m}}{m!}\sum_{l=0}^{\infty}\binom{m+\alpha+l}{l}t^{l} \nonumber \\
	&\quad =\sum_{n=0}^{\infty}\bigg(\sum_{m=0}^{n}(1)_{m,\lambda}(-1)^{m}x^{m}\frac{1}{m!}\binom{m+\alpha+n-m}{n-m}\bigg)t^{n}\nonumber \\
	&\quad =\sum_{n=0}^{\infty}\bigg(\sum_{m=0}^{n}(1)_{m,\lambda}(-1)^{m}x^{m}\frac{1}{m!}\binom{n+\alpha}{n-m}\bigg)t^{n}\nonumber .
\end{align}
Therefore, by \eqref{8} and \eqref{9}, we obtain the following theorem.
\begin{theorem}
For $n\ge 0$, we have
\begin{displaymath}
	L_{n,\lambda}^{(\alpha)}(x)=\sum_{m=0}^{n}\binom{n+\alpha}{n-m}(-1)^{m}(1)_{m,\lambda}\frac{1}{m!}x^{m}.
\end{displaymath}	
\end{theorem}
Now, by using Theorem1, we observe that
\begin{align}
&\frac{d^{n}}{dx^{n}}\bigg[x^{\alpha}e_{\lambda}\bigg(-\frac{a}{x}\bigg)\bigg]=\frac{d^{n}}{dx^{n}}\bigg[\sum_{k=0}^{\infty}(1)_{k,\lambda}(-a)^{k}\frac{1}{k!}x^{\alpha-k}\bigg]\label{10}\\
&\quad=\sum_{k=0}^{\infty}(1)_{k,\lambda}\frac{(-a)^{k}}{k!}\overbrace{(\alpha-k)(\alpha-k-1)\cdots(\alpha-k-n+1)}^{n-times}x^{\alpha-k-n}\nonumber \\
&\quad=(-1)^{n}x^{\alpha-n}\sum_{k=0}^{\infty}(1)_{k,\lambda}(-a)^{k}\frac{1}{k!}(k-\alpha)(k-\alpha+1)\cdots(k-\alpha+n-1)x^{-k}\nonumber \\
&\quad= (-1)^{n}x^{\alpha-n}n!\sum_{k=0}^{\infty}(1)_{k,\lambda}\frac{(-1)^{k}}{k!}\bigg(\frac{a}{x}\bigg)^{k}\binom{k+n-\alpha-1}{n}\nonumber\\
&\quad= (-1)^{n}x^{\alpha-n}n!\sum_{k=0}^{\infty}(1)_{k,\lambda}\frac{(-1)^{k}}{k!}\bigg(\frac{a}{x}\bigg)^{k}\sum_{l=0}^{n}\binom{n-\alpha-1}{n-l}\binom{k}{l}\nonumber \\
&\quad= (-1)^{n}x^{\alpha-n}n!\sum_{l=0}^{n}\binom{n-\alpha-1}{n-l}\sum_{k=l}^{\infty}(1)_{k,\lambda}\frac{(-1)^{k}}{k!}\bigg(\frac{a}{x}\bigg)^{k}\frac{k!}{l!(k-l)!}\nonumber \\
&\quad= (-1)^{n}x^{\alpha-n}n!\sum_{l=0}^{n}\binom{n-\alpha-1}{n-l}\sum_{k=0}^{\infty}(1)_{k+l,\lambda}(-1)^{k+l}\bigg(\frac{a}{x}\bigg)^{k+l}\frac{1}{l!k!}\nonumber \\
&\quad= (-1)^{n}x^{\alpha-n}n!\sum_{l=0}^{n}\binom{n-\alpha-1}{n-l}(-1)^{l}\bigg(\frac{a}{x}\bigg)^{l}\frac{1}{l!}(1)_{l,\lambda}\sum_{k=0}^{\infty}(1-l\lambda)_{k,\lambda}\frac{(-1)^{k}}{k!}\bigg(\frac{a}{x}\bigg)^{k}\nonumber
\end{align}
\begin{align}
&\quad= (-1)^{n}x^{\alpha-n}n!\sum_{l=0}^{n}\binom{n-\alpha-1}{n-l}(-1)^{l}(1)_{l,\lambda}\bigg(\frac{a}{x}\bigg)^{l}\frac{1}{l!}e_{\lambda}^{1-l\lambda}\bigg(-\frac{a}{x}\bigg)\nonumber \\
&\quad= (-1)^{n}x^{\alpha-n}n!e_{\lambda}\bigg(-\frac{a}{x}\bigg)\sum_{l=0}^{n}\binom{n-\alpha-1}{n-l}(-1)^{l}(1)_{l,\lambda}\bigg(\frac{a}{x-a\lambda}\bigg)^{l}\frac{1}{l!}\nonumber \\
&\quad= (-1)^{n}x^{\alpha-n}n!e_{\lambda}\bigg(-\frac{a}{x}\bigg)L_{n,\lambda}^{(-\alpha-1)}\bigg(\frac{a}{x-a\lambda}\bigg).\nonumber
\end{align}
Therefore, by \eqref{10}, we obtain the following theorem.
\begin{theorem}
For $n\ge 0$, we have
\begin{displaymath}
	\frac{d^{n}}{dx^{n}}\bigg[x^{\alpha}e_{\lambda}\bigg(-\frac{a}{x}\bigg)\bigg]= (-1)^{n}x^{\alpha-n}n!e_{\lambda}\bigg(-\frac{a}{x}\bigg)L_{n,\lambda}^{(-\alpha-1)}\bigg(\frac{a}{x-a\lambda}\bigg).
\end{displaymath}	
\end{theorem}

By using Leibniz rule and Theorem 1, we have
\begin{align}
&\frac{d^{n}}{dx^{n}}\big[e_{\lambda}(-x)x^{n+\alpha}\big]\label{12}\\
&\quad =\sum_{m=0}^{n}\binom{n}{m}\bigg[\frac{d^{m}}{dx^{m}}e_{\lambda}(-x)\bigg]\bigg[\frac{d^{n-m}}{dx^{n-m}}x^{n+\alpha}\bigg]\nonumber \\
&\quad =\sum_{m=0}^{n}\binom{n}{m}(-1)^{m}(1)_{m,\lambda}e_{\lambda}^{1-m\lambda}(-x)\cdot (n+\alpha)_{n-m}x^{n+\alpha-n+m}\nonumber \\
&\quad =n!e_{\lambda}(-x)x^{\alpha}\sum_{m=0}^{n}\binom{n+\alpha}{n-m}(-1)^{m}(1)_{m,\lambda}e_{\lambda}^{-m\lambda}(-x)x^{m}\frac{1}{m!}	\nonumber \\
&\quad =n!e_{\lambda}(-x)x^{\alpha}\sum_{m=0}^{n}\binom{n+\alpha}{n-m}(-1)^{m}(1)_{m,\lambda}\bigg(\frac{x}{1-\lambda x}\bigg)^{m}\frac{1}{m!}	\nonumber \\
&\quad=n!e_{\lambda}(-x)x^{\alpha}L_{n,\lambda}^{(\alpha)}\bigg(\frac{x}{1-\lambda x}\bigg).\nonumber
\end{align}
Thus, we obtain Rodrigues' type formula for the degenerate generalized Laguerre polynomials.
\begin{theorem}[Rodrigues' type formula]
For $n\ge 0$, we have
\begin{displaymath}
\frac{x^{-\alpha}}{n!e_{\lambda}(-x)}\frac{d^{n}}{dx^{n}}\big[e_{\lambda}(-x)x^{n+\alpha}\big]=L_{n,\lambda}^{(\alpha)}\bigg(\frac{x}{1-\lambda x}\bigg). 	\end{displaymath}
\end{theorem}
For $\alpha=-1$, from Theorem 3, we have
\begin{equation}
	\frac{x}{n!e_{\lambda}(-x)}\frac{d^{n}}{dx^{n}}\big[e_{\lambda}(-x)x^{n-1}\big]=L_{n,\lambda}^{(-1)}\bigg(\frac{x}{1-x}\bigg)	. \label{12}
\end{equation}
On the other hand, by \eqref{8}, we get
\begin{equation}
e_{\lambda}\bigg(-x\frac{t}{1-t}\bigg)=\sum_{n=0}^{\infty}L_{n,\lambda}^{(-1)}(x)t^{n}.\label{13}
\end{equation}
From \eqref{7}, we can derive the following equation.
\begin{align}
e_{\lambda}\bigg(-x\frac{t}{1-t}\bigg)\ &=\ \sum_{k=0}^{\infty}(-1)^{k}(1)_{k,\lambda}x^{k}\frac{1}{k!}\bigg(\frac{t}{1-t}\bigg)^{k} \label{14} \\
&=\ \sum_{k=0}^{\infty}(-1)^{k}(1)_{k,\lambda}x^{k}\sum_{n=k}^{\infty}L(n,k)\frac{t^{n}}{n!}\nonumber \\
&=\ \sum_{n=0}^{\infty}\bigg(\sum_{k=0}^{n}(-1)^{k}(1)_{k,\lambda}x^{k}L(n,k)\bigg)\frac{t^{n}}{n!}. \nonumber
\end{align}
Thus, by \eqref{13} and \eqref{14}, we get
\begin{equation}
L_{n,\lambda}^{(-1)}(x)=\frac{1}{n!}\sum_{k=0}^{n}(-x)^{k}(1)_{k,\lambda}L(n,k),\label{15}
\end{equation}
where $L(n,k)=\binom{n-1}{k-1}\frac{n!}{k!}$ is the Lah number. \par
Therefore, we obtain the following theorem.
\begin{theorem}
For $n\ge 0$, we have
\begin{displaymath}
L_{n,\lambda}^{(-1)}(x)=\frac{1}{n!}\sum_{k=0}^{n}(-x)^{k}(1)_{k,\lambda}L(n,k)=\sum_{k=0}^{n}(1)_{k,\lambda}(-x)^{k}\frac{1}{k!}\binom{n-1}{k-1}.
\end{displaymath}	
\end{theorem}
From Theorem 1, we note that
\begin{align}
L_{n,\lambda}^{(\alpha)}(x)\ &=\ \sum_{m=0}^{n}(1)_{m,\lambda}x^{m}(-1)^{m}\frac{1}{m!}\binom{n+\alpha}{n-m}\label{16} \\
&=\ \sum_{m=0}^{n}(1)_{m,\lambda}(-x)^{m}\frac{(n+\alpha)(n+\alpha-1)\cdots(m+\alpha+1)}{m!(n-m)!}\nonumber \\
&=\	\sum_{m=0}^{n}(1)_{m,\lambda}(-x)^{m}\frac{(n+\alpha)(n+\alpha-1)\cdots(m+\alpha+1)(m+\alpha)\cdots(\alpha+1)}{m!(n-m)!(m+\alpha)\cdots(\alpha+1)}\nonumber \\
&=\ \sum_{m=0}^{n}(1)_{m,\lambda}(-x)^{m}\frac{1}{m!(n-m)!}\frac{(n+\alpha)_{n}}{(m+\alpha)_{m}}\nonumber \\
&=\ \sum_{m=0}^{n}(1)_{m,\lambda}(-x)^{m}\frac{1}{m!(n-m)!} \frac{(n+\alpha)_{n}}{(m+\alpha)_{m}}\frac{\Gamma(\alpha+1)}{\Gamma(\alpha+1)}\nonumber \\
&=\ \Gamma(n+\alpha+1)\sum_{m=0}^{n}(1)_{m,\lambda}\frac{(-x)^{m}}{m!(n-m)!}\frac{1}{\Gamma(m+\alpha+1)}.\nonumber
\end{align}
Thus, by \eqref{16}, we get
\begin{equation}
L_{n,\lambda}^{(\alpha)}(x)=\frac{\Gamma(n+\alpha+1)}{\Gamma(n+1)}\sum_{m=0}^{n}(1)_{m,\lambda}\binom{n}{m}(-x)^{m}\frac{1}{\Gamma(m+\alpha+1)}.\label{17}
\end{equation}
In particular, $\alpha=-1$, we have
\begin{align}
L_{n,\lambda}^{(-1)}(x) &=\ \frac{1}{n}\sum_{m=0}^{n}(1)_{m,\lambda}\binom{n}{m}(-x)^{m}\frac{1}{\Gamma(m)}\label{18} \\
&=\ \sum_{m=1}^{n}(1)_{m,\lambda}(-x)^{m}\frac{1}{m!}\binom{n-1}{m-1}.\nonumber	
\end{align}
Now, we observe that
\begin{align}
\frac{d^{n}}{dx^{n}}e_{\lambda}\bigg(\frac{1}{x}\bigg)&=\frac{d^{n}}{dx^{n}}\sum_{k=0}^{\infty}(1)_{k,\lambda}\frac{1}{k!}\bigg(\frac{1}{x}\bigg)^{k}\label{19}\\
&= \sum_{k=0}^{\infty}(1)_{k,\lambda}\frac{1}{k!}(-1)^{n}\langle k\rangle_{n}x^{-n-k} \nonumber \\
&= \sum_{k=0}^{\infty}(1)_{k,\lambda}\frac{(-1)^{n}}{k!}\sum_{l=0}^{n}L(n,l)(k)_{l}x^{-n-k}\nonumber \\
&=\sum_{l=0}^{n}(-1)^{n}L(n,l)\sum_{k=0}^{\infty}(1)_{k,\lambda}\frac{(k)_{l}}{k!}x^{-n-k}\nonumber \\
&=(-1)^{n}x^{-n}\sum_{l=0}^{n}L(n,l)\sum_{k=l}^{\infty}\frac{(1)_{k,\lambda}}{k!}(k)_{l}x^{-k}\nonumber\\
& =(-1)^{n}x^{-n}\sum_{l=0}^{n}L(n,l)\sum_{k=l}^{\infty}\frac{k!}{k!(k-l)!}(1)_{k,\lambda}x^{-k} \nonumber \\
& =(-1)^{n}x^{-n}\sum_{l=0}^{n}L(n,l)x^{-l}\sum_{k=0}^{\infty}\frac{(1)_{k+l,\lambda}}{k!}x^{-k} \nonumber\\
&=(-1)^{n}x^{-n}\sum_{l=0}^{n}L(n,l)x^{-l}(1)_{l,\lambda}\sum_{k=0}^{\infty}\frac{(1-l\lambda)_{k,\lambda}}{k!}x^{-k}	\nonumber \\
& =(-1)^{n}x^{-n}\sum_{l=0}^{n}L(n,l)x^{-l}(1)_{l,\lambda}e_{\lambda}^{1-l\lambda}\bigg(\frac{1}{x}\bigg)\nonumber \\
&=(-1)^{n}x^{-n}e_{\lambda}\bigg(\frac{1}{x}\bigg)\sum_{l=0}^{n}L(n,l)(1)_{l,\lambda}x^{-l}\bigg(1+\frac{\lambda}{x}\bigg)^{-l}\nonumber \\
&=(-1)^{n}x^{-n}e_{\lambda}\bigg(\frac{1}{x}\bigg)\sum_{l=0}^{n}L(n,l)(1)_{l,\lambda}\bigg(\frac{1}{x+\lambda}\bigg)^{l}.\nonumber
\end{align}
Therefore, by \eqref{19}, we obtain the following theorem.
\begin{theorem}
For $n\ge 1$, we have
\begin{displaymath}
	\frac{d^{n}}{dx^{n}}e_{\lambda}\bigg(\frac{1}{x}\bigg)= (-1)^{n}x^{-n}e_{\lambda}\bigg(\frac{1}{x}\bigg)\sum_{l=0}^{n}L(n,l)(1)_{l,\lambda}\bigg(\frac{1}{x+\lambda}\bigg)^{l}.
\end{displaymath}
\end{theorem}
Since
\begin{displaymath}
L(n,k)=\binom{n-1}{k-1}\frac{n!}{k!}=\binom{n-1}{k-1}\frac{n!}{(n-k)!k!}(n-k)!=\binom{n-1}{k-1}\binom{n}{k}(n-k)!.	
\end{displaymath}
Thus, we have the following corollary.
\begin{corollary}
For $n\ge 1$, we have
\begin{displaymath}
\frac{d^{n}}{dx^{n}}e_{\lambda}\bigg(\frac{1}{x}\bigg)=x^{-n}(-1)^{n}e_{\lambda}\bigg(\frac{1}{x}\bigg)\sum_{l=1}
^{n}\binom{n-1}{l-1}\binom{n}{l}(1)_{l,\lambda}(n-l)!\bigg(\frac{1}{x+\lambda}\bigg)^{l}.
\end{displaymath}
\end{corollary}

\section{Degenerate Poisson random variables}
Let $X$ be the Poisson random variable with parameter $\alpha(>0)$. Then the probability mass function of $X$ is given by
\begin{displaymath}
	p(i)=P\{X=i\}=\frac{\alpha^{i}}{i!}e^{-\alpha},\quad (i=0,1,2,\dots).
\end{displaymath}
It is easy to show that
\begin{displaymath}
E\big[(X)_{n}\big]=\sum_{k=0}^{\infty}(k)_{n}p(k)=e^{-\alpha}\sum_{k=n}^{\infty}\frac{\alpha^{k}}{(k-n)!}=e^{-\alpha}\alpha^{n}\sum_{k=0}^{\infty}\frac{\alpha^{k}}{k!}=\alpha^{n}.
\end{displaymath}
Thus, we note that
\begin{displaymath}
E\bigg[\binom{X}{n}\bigg]=\frac{\alpha^{n}}{n!},\quad(n=0,1,2,\dots).	
\end{displaymath}
Let $X_{\lambda}$ be the degenerate Poisson random variable with parameter $\alpha(>0)$. Then the probability mass function of $X_{\lambda}$ is given by
\begin{displaymath}
p(i)=P\{X_{\lambda}=i\}=e_{\lambda}^{-1}(\alpha)\frac{\alpha^{i}}{i!}(1)_{i,\lambda},\quad (i=0,1,2,\dots),\quad (\mathrm{see}\ [12]).
\end{displaymath}
Then the following falling factorial moment is given by
\begin{align}
E\Big[\big(X_{\lambda}\big)_{n}\Big]&=\sum_{k=0}^{\infty}(k)_{n}p(k)=\sum_{k=0}^{\infty}(k)_{n}\frac{e_{\lambda}^{-1}(\alpha)}{k!}\alpha^{k}(1)_{k,\lambda}\label{19-1} \\
&=e_{\lambda}^{-1}(\alpha)\sum_{k=n}^{\infty}\frac{k(k-1)\cdots(k-n+1)(k-n)!}{k!(k-n)!}\alpha^{k}(1)_{k,\lambda} \nonumber\\
&=e_{\lambda}^{-1}(\alpha)\sum_{k=0}^{\infty}\alpha^{k+n}\frac{1}{k!}(1)_{k+n,\lambda}\nonumber\\
&=\alpha^{n}e_{\lambda}^{-1}(\alpha)\sum_{k=0}^{\infty}(1)_{n,\lambda}(1-n\lambda)_{k,\lambda}\frac{\alpha^k}{k!}\nonumber \\
&=\alpha^{n}e_{\lambda}^{-1}(\alpha)(1)_{n,\lambda}e_{\lambda}^{1-n\lambda}(\alpha)=\alpha^{n}(1)_{n,\lambda}\bigg(\frac{1}{1+\alpha\lambda}\bigg)^{n}\nonumber.
\end{align}
Assume that $X_{\lambda}$ is the Poisson random variable with parameter $\frac{1}{\alpha}(>0)$. Then, by using \eqref{19-1}, we obtain
\begin{align*}
\frac{d^{n}}{d\alpha^{n}}\Bigg(e_{\lambda}\bigg(\frac{1}{\alpha}\bigg)\Bigg)&=(-1)^{n}\alpha^{-n}\sum_{l=0}^{n}L(n,l)\sum_{k=0}^{\infty}\frac{(k)_{l}}{k!}(1)_{k,\lambda}\alpha^{-k}\\
&= (-1)^{n}\alpha^{-n}\sum_{l=0}^{n}L(n,l)e_{\lambda}\bigg(\frac{1}{\alpha}\bigg) e_{\lambda}^{-1}\bigg(\frac{1}{\alpha}\bigg)\sum_{k=0}^{\infty}\frac{(k)_{l}}{k!}(1)_{k,\lambda}\alpha^{-k} \\
&= (-1)^{n}\alpha^{-n}\sum_{l=0}^{n}L(n,l)e_{\lambda}\bigg(\frac{1}{\alpha}\bigg)E\Big[\big(X_{\lambda}\big)_{l}\Big] \\
&= (-1)^{n}\alpha^{-n}e_{\lambda}\bigg(\frac{1}{\alpha}\bigg)\sum_{l=0}^{n}L(n,l)(1)_{l,\lambda}\bigg(\frac{1}{\alpha+\lambda}\bigg)^{l}.
 \end{align*}

\section{Derivatives of degenerate Laguerre polynomials}
Let us consider the sequence $y_{n,\lambda}(x)$ which are given by
\begin{equation}
A(t)e_{\lambda}\bigg(-x\frac{t}{1-t}\bigg)=\sum_{n=0}^{\infty}y_{n,\lambda}(x)t^{n},\label{20}
\end{equation}
where $A(t)$ is an invertible series. \par
Note that $y_{0,\lambda}(x)=A(0)$ is a constant.  We now set $\displaystyle F_{\lambda}=F_{\lambda}(x,t)=A(t)e_{\lambda}\bigg(-\frac{x}{1-t}t\bigg)\displaystyle $.\par   From \eqref{20}, we note that
\begin{align}
\frac{\partial}{\partial x}F_{\lambda}\ &=\ A(t)\bigg(-\frac{t}{1-t}\bigg)e_{\lambda}^{1-\lambda}\bigg(-\frac{xt}{1-t}\bigg) \label{21} \\
&=\ -\bigg(\frac{t}{(1-t)-x\lambda t}\bigg)A(t)e_{\lambda}\bigg(-\frac{xt}{1-t}\bigg).\nonumber
\end{align}
By \eqref{21}, we get
\begin{equation}
\frac{\partial}{\partial x}F_{\lambda}-(1+x\lambda)t\frac{\partial}{\partial x}F_{\lambda}=-tF_{\lambda}. \label{22}	
\end{equation}
From \eqref{20} and \eqref{22}, we can derive the following equation
\begin{equation}
\sum_{n=1}^{\infty}y_{n,\lambda}^{\prime}(x)t^{n}-(1+x\lambda)\sum_{n=1}^{\infty}y_{n-1,\lambda}^{\prime}(x)t^{n}=-\sum_{n=1}^{\infty}y_{n-1,\lambda}(x)t^{n}.\label{23}
\end{equation}
By comparing the coefficients on both sides of \eqref{23}, we get
\begin{equation}
y_{n,\lambda}^{\prime}(x)-(1+x{\lambda})y_{n-1,\lambda}^{\prime}(x)=-y_{n-1,\lambda}(x),\quad (n\ge 1),\label{24}
\end{equation}
where $y_{n,\lambda}^{\prime}(x)=\frac{d}{dx}y_{n,\lambda}(x)$. \par
Now, we observe that
\begin{align}
=\frac{1}{(1-t)-x\lambda t}&=-\frac{t}{1-t}\Bigg(\frac{1}{1-\frac{x\lambda}{1-t}t}\Bigg)=-\sum_{l=0}^{\infty}x^{l}\lambda^{l}\bigg(\frac{t}{1-t}\bigg)^{l+1} \nonumber \\
&=-\sum_{l=1}^{\infty}x^{l-1}\lambda^{l-1}\bigg(\frac{t}{1-t}\bigg)^{l} \label{25} \\
&=-\sum_{l=1}^{\infty}x^{l-1}\lambda^{l-1}t^{l}\sum_{m=0}^{\infty}\binom{m+l-1}{m}t^{m}\nonumber \\
&=-\sum_{k=1}^{\infty}\bigg(\sum_{l=1}^{k}x^{l-1}\lambda^{l-1}\binom{n-1}{n-l}\bigg)t^{k}.\nonumber
\end{align}
From \eqref{21} and \eqref{25}, we can derive the following equation \eqref{26}
\begin{align}
\sum_{n=1}^{\infty}y_{n,\lambda}^{\prime}(x)t^{n}&=\frac{\partial}{\partial x}F_{\lambda}=-\bigg(\frac{t}{1-t-x\lambda t}\bigg)A(t)e_{\lambda}\bigg(-\frac{xt}{1-t}\bigg)\label{26} \\
&=\sum_{k=1}^{\infty}\bigg(-\sum_{l=1}^{k}x^{l-1}\lambda^{l-1}\binom{k-1}{k-l}\bigg)t^{k}\sum_{m=0}^{\infty}y_{m,\lambda}t^{m} \nonumber \\
&=\sum_{n=1}^{\infty}\bigg(-\sum_{k=1}^{n}\sum_{l=1}^{k}x^{l-1}\lambda^{l-1}\binom{k-1}{k-l}y_{n-k,\lambda}(x)\bigg)t^{n}.\nonumber
\end{align}
Thus, by comparing the coefficients on both sides of \eqref{26}, we get
\begin{equation}
y_{n,\lambda}^{\prime}(x)=-\sum_{k=1}^{n}\sum_{l=1}^{k}x^{l-1}\lambda^{l-1}\binom{k-1}{k-l}y_{n-k,\lambda}(x).\label{27}	
\end{equation}
Therefore, we obtain the following theorem.
\begin{theorem}
Let
\begin{displaymath}
A(t)e_{\lambda}\bigg(-\frac{xt}{1-t}\bigg)=\sum_{n=0}^{\infty}y_{n,\lambda}(x)t^{n},	
\end{displaymath}
	where $A(t)$ is an invertible series.
For $n\ge 1$, we have
\begin{displaymath}
y^{\prime}_{n,\lambda}(x)=(1+x\lambda)y_{n-1,\lambda}^{\prime}(x)-y_{n-1,\lambda}(x)
\end{displaymath}
and
\begin{displaymath}
	y^{\prime}_{n,\lambda}(x)=-\sum_{k=1}^{n}\sum_{l=1}^{k}x^{l-1}\lambda^{l-1}\binom{k-1}{k-l}y_{n-k,\lambda}(x),
\end{displaymath}
where $\displaystyle y_{n}^{\prime}(x)=\frac{d}{dx}y_{n,\lambda}(x)\displaystyle$.
\end{theorem}
Let us take $A(t)=(1-t)^{-\alpha-1}$. Then we have
\begin{equation}
\sum_{n=0}^{\infty}y_{n,\lambda}(x)t^{n}=(1-t)^{-\alpha-1}e_{\lambda}\bigg(-\frac{x}{1-x}t\bigg)=\sum_{n=0}^{\infty}L_{n,\lambda}^{(\alpha)}(x)t^{n}.\label{28}
\end{equation}
Thus, we note that $y_{n,\lambda}(x)=L_{n,\lambda}^{(\alpha)}(x),\ (n\ge 0)$. \par
Therefore, by Theorem 6 and \eqref{28}, we obtain the following corollary.
\begin{corollary}
	For $n\ge 1$, we have
	\begin{displaymath}
	\frac{d}{dx}L^{(\alpha)}_{n,\lambda}(x)=(1+x\lambda)\frac{d}{dx}L_{n-1,\lambda}^{(\alpha)}(x)-L_{n-1,\lambda}^{(\alpha)}(x)
	\end{displaymath}
	and
	\begin{displaymath}
		\frac{d}{dx}L^{(\alpha)}_{n,\lambda}(x)=-\sum_{k=1}^{n}\sum_{l=1}^{k}x^{l-1}\lambda^{l-1}\binom{k-1}{k-l}L^{(\alpha)}_{n-k,\lambda}(x).
	\end{displaymath}
\end{corollary}

\section{Conclusion}
In this paper, we introduced the degenerate generalized Laguerre polynomials, which are the first degenerate versions of the orthogonal polynomials, and derived some results related to those polynomials and Lah numbers. Some of the results are an explicit expression, Rodrigues' type formula and some expressions for the derivatives of the degenerate generalized Laguerre polynomials. \par
As one of our future projects, we would like to continue to study degenerate versions of many special numbers and polynomials.

\vspace{0.2in}

\noindent{\bf{Acknowledgments:}} Not applicable

\vspace{0.2in}

\noindent{\bf{Funding:}}
This work was supported by the Basic Science Research Program, the National Research Foundation of Korea, (NRF-2021R1F1A1050151).
\vspace{0.2in}

\noindent{\bf {Availability of data and materials:}}
No data were used to support this study.

\vspace{0.1in}

\noindent{\bf {Competing interests:}}
The authors declare no conflict of interest.

\vspace{0.1in}

\noindent{\bf{Authors' contributions:}} The authors contributed equally to the work. All
authors read and approved the final manuscript.

\end{document}